\documentclass[notitlepage,11pt]{article}
\usepackage{amsfonts,amsmath,amssymb,mathrsfs,url, 
geometry,esint,
comment
}

\catcode`\@=11
\@addtoreset{equation}{section}

\catcode`\@=12

\allowdisplaybreaks

\newtheorem{Theorem}{Theorem}[section]
\newtheorem{Lemma}[Theorem]{Lemma}
\newtheorem{Proposition}[Theorem]{Proposition}

\newtheorem{Remark}[Theorem]{Remark}
\newtheorem{Definition}[Theorem]{Definition}

\usepackage[T3,T1]{fontenc}
\DeclareSymbolFont{tipa}{T3}{cmr}{m}{n}
\DeclareMathAccent{\invbreve}{\mathalpha}{tipa}{16}

\usepackage{color}

\def\R{\mathbb{R}}

\def\f{\varphi}

\def\intl{\int\limits}
\def\fintl{\fint\limits}
\def\iintl{\iint\limits}

\def\irn{\int\limits_{\R^n}}

\def\Grho{{\Gamma(\rho)}}

\def\Dsn{\left(-\Delta_{{\rm \R^n}}\right)^{\!s}} 

\def\Dsd{\left(-\Delta_{{\rm \R^d}}\right)^{\!s}} 
\def\Dshalfd{\left(-\Delta_{{\rm \R^d}}\right)^{\!\frac{s}{2}}}

\def\DsUno{\left(-\Delta_{{\rm \R^k}}\right)^{\!s}} 

\def\Dsnp{\left(-\Delta_{{\rm \R^{n+k}}}\!\right)^{\!s}} 

\def\E{{\mathcal E}}
\def\G{{\mathscr G}}

\def\Ht{{\widetilde H}}

\def\proof{\noindent{\textbf{Proof. }}}
\def\QED{\hfill {$\square$}\goodbreak \medskip}

\begin{document}

\title 
{Asymptotic analysis of the Dirichlet fractional Laplacian in domains becoming unbounded}

\author{Vincenzo Ambrosio
\footnote{DIISM, Universit\`a Politecnica delle Marche,
via Brecce Bianche 12, 60131 Ancona, Italy. Email: {v.ambrosio@univpm.it}} ,~
Lorenzo Freddi
\footnote{DMIF, Universit\`a di Udine,
via delle Scienze 206, 33100 Udine, Italy. Email: {lorenzo.freddi@uniud.it}. Partially supported by  PRID project {\em PRIDEN}
}~
and Roberta Musina
\footnote{DMIF, Universit\`a di Udine,
via delle Scienze 206, 33100 Udine, Italy. Email: {roberta.musina@uniud.it}, 
Partially supported by PRID project {\em VAPROGE}
}~ 
}

\date{}

\maketitle

\begin{abstract}
{
{\small In this paper we analyze the asymptotic behavior of the Dirichlet fractional Laplacian $(-\Delta_{\R^{n+k}})^{s}$, with $s\in (0, 1)$, on bounded domains in $\R^{n+k}$ that become unbounded in the last $k$-directions. 
A dimension reduction phenomenon is observed and described via $\Gamma$-convergence.}}

\vskip0.5cm

\noindent
\textbf{Keywords:} {Fractional Laplacian, Variational methods, Asymptotic analysis, $\Gamma$-convergence}

\medskip\noindent
\textbf{2010 Mathematics Subject Classfication:} 35R11; 35B40; 49J45.
\end{abstract}

\section{Introduction}

The asymptotic analysis for elliptic problems on domains becoming 
unbounded in one or more directions has been extensively studied in the last decades because of its importance in applications.
In dealing for instance with elasticity properties or diffusion processes in $3$-dimensional  tubes having $2$-dimensional section 
$\omega$ and comparably long length $2\ell$, one can ask whether it is allowed to consider the same problem 
on an infinite tube with section $\omega$, in order to avoid the technical difficulties produced 
by possible "boundary effects" at the ends of the tube. The classical {\em de Saint-Venant principle}
\cite{SV} (see also \cite{Love, Toupin}) is one of the most famous mathematical models based on this basic idea (see also
the examples in \cite{Ch, Ch1, Ch3} and  references therein).
A good understanding of the asymptotic behaviour as $\ell\to \infty$ of the differential operator involved
in the model is crucially needed in order to make this approach rigorous from the mathematical point of view. 

To our knowledge, only few papers deal with the asymptotic analysis for equations involving 
non local differential operators on domains becoming infinitely long, despite
the interest towards nonlocal operators 
has considerably grown in recent years. Apart from their engaging and challenging theoretical structure,
they are largely used to model a variety of phenomena in presence of
long-range interactions. Examples arise from many fields, including 
like phase transitions, crystal dislocation, optimization, anomalous diffusion, 
semipermeable membranes and flame propagation. 

In the present paper we deal with the  behaviour of the Dirichlet 
fractional Laplacian on varying $(n+k)$-dimensional domains  
which become unbounded in the last $k$-directions. For related results we cite
\cite{ChowR, Y}, where $k=1$ is assumed, and the recent paper \cite{ChCRS}, 
where the regional (or restricted) Dirichlet  fractional Laplace operator is considered. 

More precisely, we consider the domains
$$
\Omega^{n+k}_\ell =\omega^n\times B^k_\ell\subset\R^n\times\R^k\,, \quad B^k_\ell=\ell B^k_1\,,
$$
where $n,k \ge 1$, $\omega^n$ is a given bounded and Lipschitz domain in $\R^n$, and
$B^k_1\subset \R^k$ is a bounded, open and convex neighbourhood of the origin.
For instance, $B^k_\ell$ can be the rectangle $(-\ell,\ell)^k$ or the Euclidean
ball of radius $\ell$ about the origin (from the point of view of applications,
the most interesting case is when $k=1$ and $\Omega^{n+1}_\ell$ is the cylinder
$\omega^n\times(-\ell,\ell)$ in $\R^{n+1}$). Notice that the increasing domains $\Omega^{n+k}_\ell$
give, in the limit as $\ell\to\infty$, 
the infinitely long 
$(n+k)$-dimensional cylinder
$$
\Omega^{n+k}_\infty = \omega^n\times\R^k\,.
$$

We denote by $\Dsnp$ the fractional Dirichlet Laplace operator of order $s\in(0,1)$ on
its natural domain $\widetilde H^s(\Omega^{n+k}_\ell)$ (see Section \ref{S:preliminaries} for 
notation, definitions and details), and we study its behaviour as $\ell\to\infty$. 

First, we investigate the relations between the Poincar\'e constant on $\Omega^{n+k}_\ell$, 
$$
\lambda^s(\Omega^{n+k}_\ell) =
\inf_{u\in \widetilde H^s(\Omega^{n+k}_\ell)\atop u\neq 0}
\frac{\langle\Dsnp u,u\rangle}{\|u\|^2_{L^2(\Omega^{n+k}_\ell)}}\,,\quad \ell\in(0,\infty]\,,
$$
and the Poincar\'e constants on $\omega^n$, $B^k_\ell$, namely,
$$
\lambda^s(\omega^n) = \inf_{u\in \widetilde H^s(\omega^n)\atop u\neq 0}
\frac{\langle\Dsn u,u\rangle}{\|u\|^2_{L^2(\omega^n)}}\,,\quad
\lambda^s(B^k_\ell) = \inf_{u\in \widetilde H^s(B^k_\ell)\atop u\neq 0}
\frac{\langle\DsUno u,u\rangle}{\|u\|^2_{L^2(B^k_\ell)}}\,,~~\ell\in(0,\infty)\,.
$$
A simple rescaling argument shows that $\lambda^s(B^k_\ell)=\ell^{-2s}{\lambda^s(B^k_1)}$,
see $iii)$ in Proposition \ref{P:lambdaB}.

In the local case $s=1$, derivatives along orthogonal directions do not interact and we have the splitting
$$
\begin{aligned}
\langle(-\Delta_{\R^{n+k}}) u,u\rangle=&~\intl_{\R^k}dt\irn|\nabla_x u|^2\,dx+\irn dx\intl_{\R^k}|\nabla_t u|^2\,dt\\
=&~
\intl_{\R^k}\langle (-\Delta_{\R^{n}}) u(~\!\cdot~\!,t),u(~\!\cdot~\!,t)\rangle\,dt+
\irn \langle (-\Delta_{\R^{k}}) u(x,~\!\cdot~\!),u(x,~\!\cdot~\!)\rangle\,dx
\end{aligned}
$$
for any $u\in  \widetilde H^1(\Omega^{n+k}_\ell)=H^1_0(\Omega^{n+k}_\ell)$. As a consequence, one  easily get
the equalities
$$
\lambda^1(\Omega^{n+k}_\infty)=\lambda^1(\omega^n)\,,
\qquad
\lambda^1(\Omega^{n+k}_\ell)= \lambda^1(\omega^n)+\frac{\lambda^1(B^k_1)}{\ell^2}
\,.
$$
If $s\in(0,1)$ the picture is different and technically more involved, due to the nonlocal character of the operator $\Dsnp$. In Section \ref{S:preliminaries} we first show that
\begin{multline}
\label{eq:per_tilde}
\langle\Dsnp u,u\rangle<\\
< \intl_{\R^k}\langle \Dsn u(~\!\cdot~\!,t),u(~\!\cdot~\!,t)\rangle\,dt+
\intl_{\R^n}\langle \DsUno u(x,~\!\cdot~\!),u(x,~\!\cdot~\!)\rangle\,dx
\end{multline}
for any nontrivial $u\in  \widetilde H^s(\Omega^{n+k}_\ell)$. Then we 
compute the Poincar\'e constant for the quadratic form on
the right hand side of (\ref{eq:per_tilde}) and prove our first main result, that can be stated as follows.

\begin{Theorem}
\label{T:Dirichlet}
Let $s\in(0,1)$. For any $\ell\in(0,\infty)$ we have
\begin{itemize}
\item[$i)$] $\lambda^s(\Omega^{n+k}_\infty)=\lambda^s(\omega^n)$\,,
\item[$ii)$] $\lambda^s(\omega^n)< \lambda^s(\Omega^{n+k}_\ell) < \lambda^s(\omega^n)+\dfrac{\lambda^s(B^k_1)}{\ell^{2s}}$\,.
\end{itemize}
In particular,
$$
\displaystyle{\lim_{\ell\to\infty}\lambda^s(\Omega^{n+k}_\ell)=\lambda^s(\omega^n)}~\!.
$$
\end{Theorem}

\medskip

In the second part of the paper we deal with Dirichlet problems for
the operator $\Dsnp$ on the expanding domains $\Omega^{n+k}_\ell$. Our main result
includes
\begin{equation}
\label{eq:Gell}
\begin{cases}
\Dsnp u_\ell =	f_\infty&\text{in $\Omega^{n+k}_\ell$}\\
u_\ell\equiv 0&\text{in $\R^{n+k}\setminus\overline\Omega^{n+k}_\ell$},
\end{cases}
\end{equation}
where $f_\infty=f_\infty(x)\in L^2(\Omega^n_\ell)$ is a given load which is constant in the $t$-variable.
By elementary variational arguments, problem
(\ref{eq:Gell}) admits a unique solution  $u_\ell=u_\ell(x,t)\in \widetilde H^s(\Omega^{n+k}_\ell)$.

A description of the asymptotic behaviour of $u_\ell$ in case
$k=1$, $s\in(\frac12,1)$, has been given in 
\cite[Theorem 1.2]{ChowR}. It turns out that 
$\|u_\ell- u_\infty\|_{L^2(\Omega_{\alpha\ell})}\to 0$ for any $\alpha\in(0,1)$, where $u_\infty\in\widetilde H^s(\omega^n)$ is the unique solution to
\begin{equation}
\label{eq:Ginfty}
\begin{cases}
\Dsn u_\infty =	f_\infty&\text{in $\omega^n$}\\
u_\infty\equiv 0&\text{in $\R^{n}\setminus\overline\omega^n$}.
\end{cases}
\end{equation}
It has to be noticed that the assumptions $k=1, s\in(\frac12,1)$ play a crucial role
in the argument in \cite{ChowR}, and can not be easily improved.

Once again, as well as in the local case (see \cite{Ch2}), a {\em dimension reduction} effect is observed: roughly speaking, while the section $\omega^n$ becomes somehow negligible for $\ell$  large, the
 problems (\ref{eq:Gell}) on $\Omega_\ell^{n+k}$ converge to the problem (\ref{eq:Ginfty}), which is settled on $\omega^n$. 

To perform our asymptotic analysis we adopt a variant of the notion of De Giorgi's $\Gamma$-convergence \cite{DeGiorgi}, which has been
proposed in \cite{ABP94} to study dimension reduction problems. For general definitions and 
results about $\Gamma$-convergence theory we refer to the monographs  \cite{DM} by
Dal Maso and \cite{Braides} by Braides.

A preliminary compactness analysis (see Lemma \ref{L:equicoercive}) suggests to compute the $\Gamma$-limit of the sequence of problems \eqref{eq:Gell} with respect to a convergence defined through the {\em averaging operators}
\begin{equation}
\label{eq:average}
\rho_\ell(v)(x)=\fint\limits_{B_\ell^k} v(x,t)~\!dt=\frac{1}{\ell^k|B_1^k|}\int\limits_{B_\ell^k} v(x,t)~\!dt\,,\quad v\in L^2(\Omega_\ell^{n+k})~\!.
\end{equation}
The use of convergences based on suitable averages is quite usual in the literature on dimension reduction problems, see e.g.\ \cite{ABP94, CiarletII, FP04, FP04bis}.

The next result is a particular case of the more general Theorem \ref{T:Gaeta02}.

\begin{Theorem}
\label{T:Gaeta0}
For a given $f_\infty\in  L^2(\omega^n)$, let  $u_\ell\in \widetilde H^s(\Omega^{n+k}_\ell)$,  $u_\infty\in\widetilde H^s(\omega^n)$ be the solutions to
problems (\ref{eq:Gell}), (\ref{eq:Ginfty}), respectively. Then
$$\rho_\ell(u_\ell)\to u_\infty \quad\text{strongly  in $\widetilde H^s(\omega^n)$, as $\ell\to\infty$\,.} 
$$
\end{Theorem}
Notice that 
we can handle also the case $s\in(0,\frac12]$ and obtain the strong $\widetilde H^s$ convergence of the averages.

\section{Preliminaries and proof of Theorem \ref{T:Dirichlet}}
\label{S:preliminaries}

We start by recalling some facts about the fractional Laplacian. For basic definitions and results, such as density and Rellich-type theorems,
our reference is the monograph \cite{Tr} by Triebel.

Let $s>0$ and take an integer $d\ge 1$. The $s$-Laplacian $\Dsd u$ of $u\in C^\infty_0(\R^d)$ is the smooth function 
$$
\Dsd u=\mathcal F_d^{-1}\big[|~\!\cdot~\!|^{2s}\mathcal F_d[u]\big]
$$
where 
$$\mathcal F_d[u](\xi)=\displaystyle{\frac{1}{(2\pi)^{n/2}}\int\limits_{\R^d} e^{-i z\cdot\xi}u(z)~\!dz}$$ 
is the Fourier transform of $u$ in $\R^d$. 
We will use the notation
$$
\E^s_d(u)= \langle\Dsd u,u\rangle=\intl_{\R^d}|\xi|^{2s}|\mathcal F_d[u]|^2~\!d\xi=
\intl_{\R^d}|\Dshalfd u|^2~\!dz~\!,
$$
where the last equality follows by Parseval formula. For $s=1$ we clearly have
$$
\E^1_d(u)=\intl_{\R^d}|\nabla u|^2~\!dz~\!.
$$
From now on, unless stated otherwise, we take $s\in(0,1)$. 

There are many equivalent definitions for the operator
$\Dsd$ and for its corresponding quadratic form $\mathcal E^s_d$, see \cite{K}. 
In particular, we will use the equality 
$$
\E^s_d(u)=\frac{C_{d,s}}{2}\iintl_{\R^d\times\R^d}\frac{(u(z)-u(\eta))^2}{|z-\eta|^{d+2s}}dz d\eta~,\quad 
C_{d,s}=\frac{s2^{2s}}{\Gamma(1-s)}~\!\frac{\Gamma\big(\frac{d+2s}{2}\big)}{\pi^\frac{d}{2}}~\!.
$$
 
If $u\in L^2(\R^d)$ then $\Dsd u$ is a well defined distribution on $\R^d$; the Sobolev space
$$
H^s(\R^d) = \big\{u\in L^2(\R^d)~|~ \Dshalfd u\in L^2(\R^d)~\!\big\}
$$
naturally inherits a Hilbertian structure with norm
$$
\|u\|^2 = \E^s_d(u)+\|u\|_{L^2(\R^d)}^2=\int\limits_{\R^d}\big|\Dshalfd u\big|^2~\!dx+ \int\limits_{\R^d}|u|^2 dx~\!.
$$

\medskip

Now, let $U\subset\R^{d}$ be a (possibly unbounded) Lipschitz domain. 
We regard at
$$
\widetilde H^s(U)=  \big\{~u\in H^s(\R^d)~|~ u\equiv 0~~\text{on $\R^d\setminus \overline U$}~\big\}
$$
as a closed subspace of $H^s(\R^n)$. For a given $f\in L^2(U)$, we say that $u$ is a weak solution 
of the Dirichlet problem
$$
\begin{cases} 
\Dsd u=f&\text{ in $U$}\\
u=0&\text{in $\R^d\setminus\overline U$}
\end{cases}
$$
if $u\in \widetilde H^s(U)$ and
$$\langle\Dsd u,\f\rangle=\intl_{\R^d}|\xi|^{2s}\mathcal F_d[u]\overline{\mathcal F_d[\f]}~\!d\xi=\intl_Uf\f~\!dz\quad \text{for any $\f\in \widetilde H^s(U)$.}
$$
If $U$ is bounded, then $\widetilde H^s(U)$ is compactly embedded into $L^2(U)$, the Poincar\'e 
constant 
$$
\lambda^s(U)=\inf_{u\in \widetilde H^s(U)\atop u\neq 0}
\frac{\E^s_d(u)}{\|u\|^2_{L^2(U)}}
$$
is positive, thus 
$$ 
\|u\|^2_{\widetilde H^s(U)} =\E^s_d(u)$$ 
is an equivalent Hilbertian norm on 
$\widetilde H^s(U)$. For convenience of the reader we sketch the proof of the next
standard result (for some of the statements below see also \cite[Proposition 9]{SeVa}, 
where $d>2s$ is assumed and more general nonlocal operators are considered).

\begin{Proposition}
\label{P:lambdaB}
Let $U\subset \R^d$ be a bounded and Lipschitz domain. 
\begin{itemize}
\item[$i)$] The Poincar\'e constant 
$\lambda^s(U)$ is the first eigenvalue of the
operator $\Dsd$ on $\widetilde H^s(U)$, and it is simple; up to a change of sign
its eigenfunction $e_1$ is lower semicontinuous and positive in $U$ and solves 
\begin{equation}
\label{eq:e1}
\Dsd u= \lambda^s(U)u\qquad\text{in $U$.}
\end{equation}
\item[$ii)$] If there exists a nonnegative, nontrivial function in $\widetilde H^s(U)$ such that 
$\Dsd u =\lambda u$ in $U$, then $\lambda=\lambda^{s}(U)$.
\item[$iii)$] Let $\ell>0$ and put $\ell U = \{\ell z~|~z\in U\}$. Then 
$\lambda^s(\ell U)=\ell^{-2s}{\lambda^s(U)}$.
\end{itemize}
\end{Proposition}

\proof
The existence of a minimizer $e_1$ for $\lambda^s(U)$ solving (\ref{eq:e1}) is immediate. 
To check that  $e_1$ can be assumed to be nonnegative, use the properties of the truncation $u\mapsto u^\pm$ in 
fractional Sobolev spaces; more precisely, it holds that 
$\mathcal E^s_d(|u|)< \mathcal E^s_d(u)$
for any sign-changing function $u\in \widetilde H^s(U)$, see for instance \cite[Theorem 6]{MN}. In particular, 
$\Dsd e_1\ge 0$ in the distributional sense on $U$.
By the strong maximum principle in \cite[Corollary 4.2]{MaxPri} (see also 
\cite[Section 2]{Sil}, where $d>2s$ is needed), the eigenfunction $e_1$ is lower semicontinuous and positive in $U$. 

To prove that $\lambda^s(U)$ is a simple eigenvalue for $\Dsd$ we argue in a standard way. We take a nontrivial solution $u$ to (\ref{eq:e1})
and notice that the function 
$$
\tilde u=u\intl_U e_1~\!dz-e_1\intl_U u~\!dz
$$
solves (\ref{eq:e1}), thus it vanishes on $U$ (otherwise, it would be a sign-changing extremal for the Poincar\'e constant,
which is impossible). Thus $u$ is proportional to $e_1$, as claimed.

To prove $ii)$ observe, as usual, that the functions $u$ and $e_1$ can not be orthogonal in $L^2(U)$ because $ue_1$ has constant sign.

It remains to prove $iii)$. Take $\ell,\ell'>0$ and any nontrivial function $v_{\ell'}\in \widetilde H^s(\ell' U)$. 
Use the function
$v_\ell(z)=v_{\ell'}\big(\frac{\ell' }{\ell}z\big)\in \widetilde H^s(\ell U)$ to estimate $\lambda^s(\ell U)$. 
Since 
$$
\E^s_d(v_\ell)=\Big(\frac{\ell}{\ell'}\Big)^{d-2s}\E^s_d(v_{\ell'})~,\quad \|v_\ell\|^2_{L^2(\ell U)}=\Big(\frac{\ell}{\ell'}\Big)^{d} \|{v_{\ell'}}\|^2_{L^2(\ell' U)}
$$
we plainly get
$$
\lambda^s(\ell U)\le\frac{\displaystyle\E^s_d(v_\ell)}{\|v_\ell\|^2_{L^2(\ell U)}}=
\Big(\frac{\ell'}{\ell}\Big)^{2s}~\!\frac{\displaystyle\E^s_d(v_{\ell'})}{\|{v_{\ell'}}\|^2_{L^2(\ell' U)}}~\!.
$$
Thus $\ell^{2s}\lambda^s(\ell U)\le (\ell')^{2s}\lambda^s(\ell' U)$, because $v_{\ell'}$ was arbitrarily chosen
in $\widetilde H^{s}(\ell' U)$. Since $\ell,\ell'$ were arbitrarily chosen as well, the conclusion follows.
\QED

\medskip

Now we take integers $n,k\ge 1$ and work in the space $H^s(\R^{n+k})$. 
We denote points in $\R^{n+k}\equiv \R^n\times\R^k$ as pairs $(x,t)$ with $x\in\R^n$ and $t\in\R^k$,
and introduce the quadratic form
\begin{equation}
\label{eq:tildeE}
\widetilde{\mathcal E}^s_{n+k}(u) =  \intl_{\R^k}\E^s_{n}(u(~\!\cdot~\!,t))~\!dt+\intl_{\R^n}\E^s_{k}(u(x,~\!\cdot~\!))~\!dx~\!.
\end{equation}
The next lemma will be crucially used in the proofs of our main results.

\begin{Lemma}
\label{L:dim}
Let $s\in(0,1)$ and $u\in H^s(\R^{n+k})$ be given. 
\begin{itemize}
\item[$i)$] $u(~\!\cdot~\!,t)\in H^{s}(\R^{k})$ for a.e. $t\in \R^k$ and $u(x,~\!\cdot~\!)\in H^s(\R^n)$ 
for a.e. $x\in \R^n$;
\item[$ii)$] If $u$ is nontrivial then
\begin{gather}
\label{eq:first}
\intl_{\R^k}\E^s_{n}(u(~\!\cdot~\!,t))\,dt< ~\E^s_{n+k}(u)\,,\quad
\intl_{\R^n}\E^s_{k}(u(x,~\!\cdot~\!))\,dx< ~\E^s_{n+k}(u)\,,\\
\label{eq:second_bis}
2^{s-1} \widetilde{\mathcal E}^s_{n+k}(u)< \E^s_{n+k}(u)< \widetilde{\mathcal E}^s_{n+k}(u).
\end{gather}
\end{itemize}
\end{Lemma}

\proof
Notice that $\mathcal F_{n+k}[u](\xi,\tau)= \mathcal F_k[\hat u^\xi](\tau)$, where
$\hat u^\xi(t)=\mathcal F_n[u(~\!\cdot~\!,t)](\xi)$. Thus, Parseval formula and Fubini's theorem give
$$
\begin{aligned}
\irn |\xi|^{2s}~\!d\xi \intl_{\R^k} \big|\mathcal F_{n+k}[u](\xi,\tau)\big|^{2}\,d\tau=&
\irn |\xi|^{2s}d\xi\intl_{\R^k}\big|\mathcal F_k[\hat u^\xi](\tau)\big|^2\,d\tau
=
\irn |\xi|^{2s}d\xi\intl_{\R^k}|\hat u^\xi(t)|^2\,dt\\
=&
\intl_{\R^k}~\!dt \irn |\xi|^{2s}\big|\mathcal F_n[u(~\!\cdot~\!,t)](\xi)\big|^2\,d\xi.
\end{aligned}
$$
We infer that
\begin{equation}
\label{eq:noCS1}
\irn |\xi|^{2s}~\!d\xi \intl_{\R^k} \big|\mathcal F_{n+k}[u](\xi,\tau)\big|^{2}\,d\tau=
\intl_{\R^k} \E^s_n(u(~\!\cdot~\!,t))\,dt.
\end{equation}
In a similar way one can check that 
\begin{equation}
\label{eq:noCS2}
\intl_{\R^k}|\tau|^{2s}~\!d\tau \irn\big|\mathcal F_{n+k}[u](\xi,\tau)\big|^{2}\,d\xi=\intl_{\R^n}\E^s_{k}(u(x,~\!\cdot~\!))\,dx\,.
\end{equation}
Since 
$$
\E^s_{n+k}(u)=\iintl_{\R^n\times\R^k}\big(|\xi|^2+|\tau|^2\big)^{s}\big|\mathcal F_{n+k}[u]\big|^{2}\,d\xi $$
the conclusion follows from (\ref{eq:noCS1}), (\ref{eq:noCS2}) and the elementary inequalities
$$
|\xi|^{2s}, |\tau|^{2s}< \big(|\xi|^2+|\tau|^2\big)^{s}\,,\quad 
2^{s-1}(|\xi|^{2s}+|\tau|^{2s})< \big(|\xi|^2+|\tau|^2\big)^{s}< |\xi|^{2s}+|\tau|^{2s}\,,
$$
that hold for every $(\xi,\tau)\in\R^{n+k}$ such that $\xi\ne0\ne\tau$ and $|\xi|\ne|\tau|$.
\QED

\begin{Remark}
If $s>1$, we have the opposite inequalities
$
2^{s-1}\widetilde{\mathcal E}^s_{n+k}(u)> 
\E^s_{n+k}(u)> \widetilde{\mathcal E}^s_{n+k}(u)$ for any nontrivial $u\in H^s(\R^{n+k})$.

\end{Remark}

We conclude this preliminary section by pointing out a consequence of Lemma \ref{L:dim}.

\begin{Lemma}
\label{L:average}
Let $s\in(0,1)$ and let $\rho_\ell : L^2(\Omega_\ell^{n+k})\to L^2(\omega^n)$ be the averaging map in (\ref{eq:average}). Then
$\rho_\ell :  \widetilde H^s(\Omega_\ell^{n+k})\to  \widetilde H^s(\omega^{n})$, and
$$
\mathcal E^s_n(\rho_\ell(v))\le \fint\limits_{B_\ell^k} \mathcal E^s_n(v(~\!\cdot~\!,t))~\!dt\le \frac{1}{|B_\ell^k|}
\mathcal E^s_{n+k}(v)\quad\text{for any $v\in \widetilde H^s(\Omega_\ell^{n+k})$.}
$$
\end{Lemma}

\proof
Take $v\in \widetilde H^s(\Omega_\ell^{n+k})$ and use the Cauchy-Schwarz inequality to get
$$
(\rho_\ell(v)(x)-\rho_\ell(v)(y))^2\le \fint\limits_{B_\ell^k}(v(x,t)-v(y,t))^2\,dt\,.
$$
Thus, Fubini's theorem gives
$$
\begin{aligned}
\mathcal E^s_n(\rho_\ell(v))=&~\frac{C_{n,s}}{2}\iintl_{\R^n\times\R^n}\frac{(\rho_\ell(v)(x)-\rho_\ell(v)(y))^2}{|x-y|^{n+2s}}dxdy\\
\le&~ \frac{C_{n,s}}{2}\fint\limits_{B_\ell^k}dt\iintl_{\R^n\times\R^n}\frac{(v(x,t)-v(y,t))^2}{|x-y|^{n+2s}}dxdy=
\fint\limits_{B_\ell^k} \mathcal E^s_n(\rho_\ell(v(~\!\cdot~\!,t))\,dt.
\end{aligned}
$$
Finally, we apply Lemma \ref{L:dim} to conclude the proof.
\QED

We now collect few additional remarks about the quadratic form $\widetilde{\mathcal E}^s_{n+k}$ in (\ref{eq:tildeE}). 

First of all, notice that $\|u\|^2=\widetilde{\mathcal E}^s_{n+k}(u)$ defines an 
equivalent Hilbertian norm 
in $\widetilde H^s(\Omega_\infty^{n+k})$ by (\ref{eq:second_bis}). Next, we study the Poincar\'e constants
$$
\widetilde{\lambda}^s(\Omega^{n+k}_\ell) = \inf_{u\in \widetilde H^s(\Omega^{n+k}_\ell)\atop u\neq 0}
\frac{\widetilde{\mathcal E}^s_{n+k}(u)}{\|u\|^2_{L^2(\Omega^{n+k}_\ell)}}~\!, \quad \ell\in(0,\infty].
$$

\begin{Lemma}
\label{L:tildeE}
It holds that  
$$
\widetilde{\lambda}^s(\Omega^{n+k}_\ell)=\lambda^s(\omega^n)+\dfrac{\lambda^s(B^k_1)}{\ell^{2s}}
\quad\text{if $\ell\in(0,\infty)$,}\quad \widetilde{\lambda}^s(\Omega^{n+k}_\infty)=\lambda^s(\omega^n)~\!.
$$
\end{Lemma}

\proof
Let $\ell\in(0,\infty)$ and take any $u\in \widetilde H^s(\Omega^{n+k}_\ell)$. By 
$i)$ in Lemma \ref{L:dim} we get
\begin{equation}
\label{eq:elementare}
\lambda^s(\omega^n)\intl_{\omega^n}|u(x,t)|^2~\!dx\le \E^s_{n}(u(~\!\cdot~\!,t))~,\quad 
\lambda^s(B_\ell^k)\intl_{B^k_\ell}|u(x,t)|^2~\!dt\le \E^s_{k}(u(x,~\!\cdot~\!))
\end{equation}
for almost every $x\in\R^n$, $t\in\R^k$. Thus, using also $iii)$ in Proposition \ref{P:lambdaB}, we obtain
\begin{eqnarray*}
&&\Bigl(\lambda^s(\omega^n)+\dfrac{\lambda^s(B^k_1)}{\ell^{2s}}\Big)\intl_{\Omega_\ell^{n+k}}|u|^2~\!dxdt=\\
&&\hspace{60pt}=
\lambda^s(\omega^n)\intl_{\R^k}dt\intl_{\omega^n}|u(x,t)|^2~\!dx+
\lambda^s(B_\ell^k)\intl_{\R^n}dx\intl_{B^k_\ell}|u(x,t)|^2~\!dt\\
&&\hspace{60pt}
\le\intl_{\R^k}\E^s_{n}(u(~\!\cdot~\!,t))~\!dt+\intl_{\R^n}\E^s_{k}(u(x,~\!\cdot~\!))~\!dx=
\widetilde{\mathcal E}^s_{n+k}(u)~\!,
\end{eqnarray*}
that implies $\lambda^s(\omega^n)+\dfrac{\lambda^s(B^k_1)}{\ell^{2s}}\le 
\widetilde{\lambda}^s(\Omega^{n+k}_\ell)$, since $u$ was arbitrarily chosen.

To prove the opposite inequality we take an eigenfunction $V \in\widetilde H^s(\omega^n)$ for the operator $\Dsn$, an
eigenfunction $v_\ell \in\widetilde H^s(B^k_\ell)$ for the operator $\DsUno$ and
test $\widetilde{\lambda}^s(\Omega^{n+k}_\ell)$ with the function  $(Vv_\ell)(x,t)=V(x)v_\ell(t)$. Since
$$
\begin{aligned}
\widetilde{\mathcal E}^s_{n+k}(Vv_\ell)=&~
\E^s_{n}(V)\|v_\ell\|^2_{L^2(B^k_\ell)} +\E^s_{k}(v_\ell)\|V\|^2_{L^2(\omega^n)}\\
=&~
\big(\lambda^s(\omega^n)+\lambda^s(B^k_\ell)\big)\|v_\ell\|^2_{L^2(B^k_\ell)}\|V\|^2_{L^2(\omega^n)}
=\big(\lambda^s(\omega^n)+\lambda^s(B^k_\ell)\big)\|Vv_\ell\|^2_{L^2(\Omega_\ell^{n+k})}~\!,
\end{aligned}
$$
we infer that $\lambda^s(\omega^n)+{\ell^{-2s}}{\lambda^s(B^k_1)}= 
\widetilde{\lambda}^s(\Omega^{n+k}_\ell)$ and that $Vv_\ell$ achieves $\widetilde{\lambda}^s(\Omega^{n+k}_\ell)$.

\medskip

It remains to prove that $\widetilde{\lambda}^s(\Omega^{n+k}_\infty)=\lambda^s(\omega^n)$. 
Since 
$\widetilde H^s(\Omega_\ell^{n+k})\subset \widetilde H^s(\Omega_\infty^{n+k})$ then
$$
\widetilde\lambda^s(\Omega_\infty^{n+k})\le \lim_{\ell\to \infty} \widetilde\lambda^s(\Omega_\ell^{n+k})=\lambda^s(\omega^n)\,.
$$
Next we integrate on $\R^k$ the first equality in (\ref{eq:elementare}) and then we use (\ref{eq:first}) and (\ref{eq:second_bis}) to get
\begin{equation}
\label{eq:gia_fatto}
\lambda^s(\omega^n)\|u\|_{L^2(\Omega_\infty^{n+k})}^2\le \E^s_{n+k}(u)
\le \widetilde \E^s_{n+k}(u)\qquad\text{for any $u\in \widetilde H^s(\Omega_\infty^{n+k})$,}
\end{equation}
which, in particular, gives $\widetilde\lambda^s(\Omega_\infty^{n+k})\ge \lambda^s(\omega^n)$.
The lemma is completely proved.
\QED

\medskip

\noindent
{\bf Proof of Theorem \ref{T:Dirichlet}}\\
We start by noticing that
(\ref{eq:gia_fatto}) gives 
$$
\lambda^s(\omega^n)\le \lambda^s(\Omega^{n+k}_\infty)~\!.
$$
The function $\ell\mapsto \lambda^s(\Omega^{n+k}_\ell)$ is decreasing, thus
$\lambda^s(\Omega^{n+k}_\infty)\le \lambda^s(\Omega^{n+k}_\ell)$. Actually
$$
~~\lambda^s(\Omega^{n+k}_\infty)< \lambda^s(\Omega^{n+k}_\ell).
$$
In fact, recall that $\lambda^s(\Omega^{n+k}_\ell)$ is attained
by some function $u_\ell \in \widetilde H^s(\Omega_\ell^{n+k})$, and notice that  $u_\ell$ can not achieve 
$\lambda^s(\Omega^{n+k}_\infty)$ because of the maximum principle in \cite{MaxPri}.

Next, we notice that 
$$\lambda^s(\Omega^{n+k}_\ell)\le  \lambda^s(\omega^n)+\dfrac{\lambda^s(B^k_1)}{\ell^{2s}}$$
because ${\mathcal E}^s_{n+k}(u)\le \widetilde{\mathcal E}^s_{n+k}(u)$ for any $u\in \widetilde H^s(\Omega^{n+k}_\ell)$ by 
Lemma \ref{L:dim}, and thanks to  Lemma \ref{L:tildeE}. In fact the strict inequality
$$
\lambda^s(\Omega^{n+k}_\ell)<  \lambda^s(\omega^n)+\dfrac{\lambda^s(B^k_1)}{\ell^{2s}}
$$
holds. Here is the simple argument. Recall that the function  in the proof of Lemma \ref{L:tildeE} achieves 
$\widetilde{\lambda}^s(\Omega^{n+k}_\ell)$. Then notice that $Vv_\ell$ it can not achieve ${\lambda}^s(\Omega^{n+k}_\ell)$,
because of the strict inequality in (\ref{eq:second_bis}). 

In conclusion, we have  the chain of inequalities
$$
\lambda^s(\omega^n)\le \lambda^s(\Omega^{n+k}_\infty)
< \lambda^s(\Omega^{n+k}_\ell)<  \lambda^s(\omega^n)+\dfrac{\lambda^s(B^k_1)}{\ell^{2s}}~\!,
$$
which implies both $i)$ and $ii)$. 
\QED

\begin{Remark}
The Poincar\'e constant for $\Omega^{n+k}_\infty$ is clearly not achieved; it is enough to use the equality $\lambda^s(\Omega^{n+k}_\infty)=\lambda^s(\omega^n)$ and the first inequality in (\ref{eq:first}).
\end{Remark}

\section{The $\Gamma$-convergence result}
\label{S:Gamma}

\def\Iscr{{{\mathscr E}}}
\def\Fscr{{{\mathscr F}}}
\def\Gscr{{{\mathscr G}}}
\def\wto{{\rightharpoonup}}

For any $\ell>0$, let $f_\ell\in L^2(\Omega^{n+k}_\ell)$ be given and let $u_\ell=u_\ell(x,t)\in \widetilde H^s(\Omega^{n+k}_\ell)$
be the unique solution to the Dirichlet problem
\begin{equation}
\label{eq:Gell2}
\begin{cases}
\Dsnp u_\ell=	f_\ell(x,t)&\text{in $\Omega^{n+k}_\ell$}\\
u_\ell\equiv 0&\text{in $\R^{n+k}\setminus\overline\Omega^{n+k}_\ell$.}
\end{cases}
\end{equation}
In this section we study the asymptotic behaviour of the solution $u_\ell$ to (\ref{eq:Gell2}) when the
loads $f_\ell$ suitably converge to a given function
$$
f_\infty=f_\infty(x)~,\qquad f_\infty \in L^2(\omega^n)~\!.
$$
More precisely, we assume that 
\begin{equation}
\label{eq:ip1}
\invbreve{f}_\ell\to f_\infty\quad \text{strongly in $L^2(\Omega_1^{n+k})$, ~~where~~ $\invbreve{f}_\ell(x,t)=f_\ell(x,\ell t)$}
\end{equation}
(as usual we regard at $f_\infty$ as a function in $L^2(\Omega_1^{n+k})$ that is constant 
in the $t$ variable).

One of the goals of the present section is to prove the next result, which includes Theorem \ref{T:Gaeta0} as a corollary.

\begin{Theorem}
\label{T:Gaeta02}
Let $\rho_\ell:\widetilde H^s(\Omega^{n+k}_\ell)\to  \widetilde H^s(\omega^n)$ be the averaging operator in (\ref{eq:average}). If (\ref{eq:ip1})
holds, then $\rho_\ell(u_\ell)\to u_\infty$  strongly in $\widetilde H^s(\omega^n)$, where $u_\infty$ solves (\ref{eq:Ginfty}).
\end{Theorem}

The argument we use to prove Theorem \ref{T:Gaeta02}
 underlines the variational nature of the $G$-type convergence result in Theorem \ref{T:Gaeta02},
 see Theorem \ref{T:Gamma_omega} below. 

We exploit the variational characterization of 
$u_\ell, u_\infty$ as the unique solutions of the convex minimization problems
$$
M_\ell=\inf_{u\in\Ht^s(\Omega^{n+k}_\ell)} \dfrac12\E^s_{n+k}(u)-\mathcal G_\ell(u)~\!,~\quad
M_\infty=\inf_{u\in\Ht^s(\omega^n)} \dfrac12\E^s_{n}(u)-\mathcal G_\infty(u)
$$
respectively, where the loading functionals $\mathcal G_\ell: \Ht^s(\Omega^{n+k}_\ell)\to \R$, $ \mathcal G_\infty:\Ht^s(\omega^n)\to\R$
are defined by
$$
\mathcal G_\ell(u)=\int\limits_{\Omega^{n+k}_\ell}f_\ell(x,t) u\,dxdt~,\quad
\mathcal G_\infty(u)= \int\limits_{\omega^n}f_\infty(x) u\,dx~\!.
$$

Since we aim to compute a variational limit of the sequence of energies, we prefer to work with functionals defined on the same domain. To do this,    
it is convenient to rescale the functionals involved in the definition of 
 $M_\ell$ as follows. 
 First of all, we introduce the transform $v\mapsto \breve v$
 $$
\breve{v}_\ell(x,t)= v\big(x,\frac{t}{\ell}\big)\,,\qquad (x,t)\in \Omega_1^{n+k}=\omega^n\times B_1^k\,,
$$
and the  functionals $\Iscr_\ell,\, \G_{\ell}: \Ht^s(\Omega_1^{n+k})\to\R$ defined by
$$
\Iscr_\ell(v)=\frac{1}{|B^k_\ell|}~\!\E^s_{n+k}(\breve{v}_\ell)~,\quad 
\G_{\ell}(v)=\frac{1}{|B^k_\ell|} {{\mathcal G_{\ell}(\breve{v}_\ell)}}\,.
$$
We explicitly have
\begin{eqnarray}
\nonumber
\Iscr_\ell(v)\!\!&=&\!\!\ell^k~\!\frac{C_{n+k,s}}{2|B^k_1|}\iintl_{\R^n\times\R^n}dxdy\!\!\iintl_{\R^k\times\R^k}
\frac{(v(x,t)-v(y,\tau))^2}{(|x-y|^2+\ell^2|t-\tau|^2)^{\frac{n+k+2s}{2}}}~\!dtd\tau\,,\\
\G_{\ell}(v)\!\!&=&\!\!\frac{1}{|B^k_1|}\intl_{\Omega_1^{n+k}} \invbreve{f}_\ell(x,t)v~\!dxdt\,.
\label{eq:load}
\end{eqnarray}
Notice that $\Iscr_\ell(v)^{1/2}$ is an equivalent Hilbertian norm on $\widetilde H^s(\Omega_1^{n+k})$. Moreover,
\begin{equation}
\label{eq:I_poincare}
\inf_{v\in \widetilde H^s(\Omega^{n+k}_1)\atop v\neq 0}
\frac{|B_1^k|\Iscr_\ell(v)}{\|v\|^2_{L^2(\Omega_1^{n+k})}}=
\inf_{u\in \widetilde H^s(\Omega^{n+k}_\ell)\atop u\neq 0}
\frac{\E^s_{n+k}(u)}{\|u\|^2_{L^2(\Omega_\ell^{n+k})}}=
\lambda^s(\Omega_\ell^{n+k})>\lambda^s(\omega^n)
\end{equation}
by Theorem \ref{T:Dirichlet}. 
Evidently, the function
$$
\invbreve{u}_\ell(x,t)= u_\ell(x,\ell t)
$$
is the unique solution to the minimization problem
$$
\invbreve{M}_\ell=\inf_{v\in \Ht^s(\Omega_1^{n+k})} I_\ell(v)\,,\quad
I_\ell(v)=\frac12\Iscr_\ell(v)-\G_\ell(v)\,\!,
$$
and we are led to investigate the ``variational convergence'' of the functionals $I_\ell$  towards the limit functional
$$
I_\infty(u) = \frac12 \mathcal E^s_n(u)-{\mathcal G}_\infty(u)\,,\quad I_\infty:\Ht^s(\omega^n)\to\R\,.
$$
We use here a variant of De Giorgi's notion of $\Gamma$-limit, that has been introduced in \cite{ABP94} to 
study dimension reduction problems and that allows to handle limiting functionals whose domain is different from that one of the approaching sequence of functionals.

The effectiveness of $\Gamma$-convergence relies on the variational property (see for instance Braides \cite{Braides} or Dal Maso \cite{DM}) that holds when the $\Gamma$-limit is computed with respect to a convergence that ensure compactness to the sequences of competing functions having bounded energy (sequential equicoercivity).
It can be seen, for instance, that in the class of problems under consideration the standard weak convergence is not suitable because of a lack of coercitivity (or compactness) (see Section \ref{sec:rem} for further remarks on this issue). 
We stress the fact that
exactly the same difficulty raises in the local case $s=1$.

This is the main reason that led us to introduce the averaging operators $\rho_\ell$.

Having in mind the definition \eqref{eq:average} of the operators $\rho_\ell :\widetilde H^s(\Omega_\ell^{n+k})\to \widetilde H^s(\omega^{n})$, we denote by 
$\rho=\rho_1$ the averaging operator on the unit ball of $\R^k$. Thus
$$
\rho(v)(x)=\fint\limits_{B_1^k} v(x,t)\,dt\,,\quad \rho : \Ht^s(\Omega_1^{n+k})\to \widetilde H^s(\omega^n)\,;
$$
compare with Lemma \ref{L:average}.
Let us point out two inequalities that will be useful in proving the subsequent compactness lemma. 
Let  $w\in \Ht^s(\Omega_1^{n+k})$ and  $\ell>0$ be given. Since  $\rho_\ell(\breve{w}_\ell)=\rho(w)$, then Lemma \ref{L:average}
(with $\breve{w}_\ell$ instead of $v$) gives
\begin{equation}
\label{eq:I_ine}
\mathcal E^s_n(\rho(w))\le\fint\limits_{B^k_1}\mathcal E^s_n(w(~\!\cdot~\!,t))~\!dt\le  \Iscr_\ell(w)\,,
\end{equation}
while, on the other hand, using \eqref{eq:I_poincare} we know that
\begin{equation}
\label{eq:I_ine2}
{{\frac{1}{|B_{1}^{k}|}}}\lambda^s({\omega^{n}})\|w\|^2_{L^2(\Omega^{n+k}_1)}\le  \Iscr_\ell(w).
\end{equation}
Both inequalities hold for any $w\in \Ht^s(\Omega_1^{n+k})$ and $\ell>0$.

We are now in position to state and proof the compactness result.

\begin{Lemma}[Equicoercivity]
\label{L:equicoercive}
Let  $v_\ell\in \Ht^s(\Omega_1^{n+k})$ be given, and assume that $I_\ell(v_\ell)$ is bounded. Then there exist a subsequence $\ell_h\to\infty$,
and a function $v\in L^2(\Omega^{n+k}_1)$, such that $\rho(v)\in \widetilde H^s(\omega^n)$ and
\begin{gather*}
v_{\ell_h}\to v\quad\text{weakly in $L^2(\Omega^{n+k}_1)$},\\
 \rho(v_{\ell_h})\to \rho(v)\quad\text{weakly in $\widetilde H^s(\omega^n)$.}
\end{gather*}
\end{Lemma}

\proof Since the sequence $\invbreve{f}_\ell$ is bounded in $L^2(\Omega_1^{n+k})$,
the Cauchy-Schwarz and Poincar\'e inequalities give
$$
I_\ell(v_\ell)\ge {{\frac{1}{2}}} \Iscr_\ell(v_\ell)-  \frac{1}{|B^k_1|}\|\invbreve{f}_\ell\|_{L^2(\Omega_1^{n+k})}\|v_\ell\|_{L^2(\Omega_1^{n+k})}
\ge {{\frac{1}{2}}} \Iscr_\ell(v_\ell)- c~\!\Iscr_\ell(v_\ell)^{1/2},
$$
where the constant $c$ does not depend on $\ell$. We infer that the sequence $\Iscr_\ell(v_\ell)$ is bounded as well.
Therefore, thanks to (\ref{eq:I_ine2}) we can find a subsequence $\ell_h$ and a function $v\in L^2(\Omega^{n+k}_1)$
such that $v_{\ell_h}\to v$ weakly in $L^2(\Omega^{n+k}_1)$.
The linear transform $\rho:L^2(\Omega^{n+k}_1)\to L^2(\omega^n)$ is continuous, hence $\rho(v_{\ell_h})\to \rho(v_0)$ weakly in $L^2(\omega^n)$. 
On the other hand, from (\ref{eq:I_ine}) we see that the sequence $\rho(v_\ell)$ is bounded in  $\widetilde H^s(\omega^n)$. Thus we can assume that 
$\rho(v_{\ell_h})\to u$ weakly in $\widetilde H^s(\omega^n)$ for some $u\in \widetilde H^s(\omega^n)$. But then, Rellich theorem guarantees
that $\rho(v_{\ell_h})\to u$ strongly in $L^2(\omega^n)$, that implies $u=\rho(v)$ and concludes the proof.
\QED

The main result of this section follows. With a slight  abuse of notation, we denote by 
$\ell\to \infty$ any given divergent sequence.

\begin{Theorem}
\label{T:Gamma_omega}
The functional $I_\infty:\widetilde H^s(\omega^n)\to \overline\R$ is the $\Grho$-limit of the
sequence $I_\ell$ on $\widetilde H^s(\omega^n)$, and write  $I_\infty=\Grho\lim\limits_{\ell\to \infty} I_\ell$,
in the sense that the following facts hold:
\begin{itemize}
\item[$i)$] for every $u\in  \widetilde H^s(\omega^n)$ and every sequence $v_\ell\in \Ht^s(\Omega_1^{n+k})$ such that 
$\rho({v}_\ell)\to u$ weakly in $\widetilde H^s(\omega^n)$, we have the ``liminf inequality''
$$
I_\infty(u)\le\liminf_{\ell\to\infty}I_\ell(v_\ell);
$$
\item[$ii)$] for every $u\in \widetilde H^s(\omega^n)$ there exists a  ``recovery sequence'' $\overline{v}_\ell\in \Ht^s(\Omega_1^{n+k})$, such that
$\rho(\overline{v}_\ell)\to u$ weakly in $\widetilde H^s(\omega^n)$ and such that the following ``limsup inequality'' holds,
$$
I_\infty(u)\ge\limsup_{h\to\infty}I_\ell(\overline{v}_\ell)~\!.
$$ 
\end{itemize}  
\end{Theorem}

\proof
Take $u\in  \widetilde H^s(\omega^n)$ and any sequence $v_\ell\in \Ht^s(\Omega_1^{n+k})$ such that 
$\rho({v}_\ell)\to u$ weakly in $\widetilde H^s(\omega^n)$. 
If $\liminf\limits_{\ell\to\infty}I_\ell(v_\ell)=+\infty$, then we are done. Otherwise, 
can assume that the sequence $I_\ell (v_\ell)$ converges. Then Lemma \ref{L:equicoercive} and Rellich theorem guarantee the existence 
of a subsequence $\ell\to\infty$ and of a function $v\in L^2(\Omega^{n+k}_1)$ such that 
$$
v_{\ell}\to v\quad\text{weakly in $L^2(\Omega^{n+k}_1)$, and $\rho(v_{\ell})\to \rho(v)=u$ in $L^2(\omega^n)$.}
$$
Therefore 
\begin{equation}
\label{eq:G}
\begin{aligned}
\G_\ell(v_\ell)=&~\intl_{\omega^n}\fintl_{B^k_1} \invbreve{f}_\ell(x,t)v_\ell~\!dxdt=
\intl_{\omega^n}\fintl_{B^k_1} {f}_\infty(x)v~\!dxdt+o(1)\\
=&~ \intl_{\omega^n}{f}_\infty(x)~\!dx\fintl_{B^k_1} v~\!dt+o(1)
=\intl_{\omega^n}{f}_\infty(x)\rho(v)~\!dx+o(1)\\
=&~\intl_{\omega^n}{f}_\infty(x)u~\!dx+o(1)=\mathcal G_\infty(u)+o(1)~\!.
\end{aligned}
\end{equation}
The weak lower semicontinuity of $\mathcal E^s_n$ on 
$\widetilde H^s(\omega^n)$ and (\ref{eq:I_ine}) give
$$
\mathcal E^s_n(u)\le \liminf_{\ell\to\infty} \mathcal E^s_n(\rho(v_{\ell}))\le \liminf_{\ell\to\infty} \Iscr_\ell (v_{\ell}),
$$
hence
$I_\infty(u)= {{\frac{1}{2}}}\mathcal E^s_n(u)-\mathcal G_\infty(u)\le \liminf\limits_{\ell\to\infty} ({{\frac{1}{2}}} \Iscr_\ell (v_{\ell})-\G_\ell(v_\ell))
=\lim\limits_{\ell\to\infty} I_\ell(v_{\ell})$, 
that concludes the proof of the liminf inequality in $i)$. 

\medskip

Next we deal with the limsup inequality. Fix a function $u\in \Ht^s(\omega^n)$. 
To construct a recovery sequence for $u$ we 
take a nonnegative cut-off function $\phi\in C^\infty_0(\R)$ with support in $(-1,1)$ and such that $\phi(0)=1$. Then we consider the sequence of functions $\f_\ell\in C^\infty_0(B^k_1)$ defined by
$$
\f_\ell(t)=\phi(|t|^{\ell}).
$$
Lebesgue's Theorem readily implies that $\f_\ell$ converges to the constant function $1$ in $L^p(B^k_1)$, for any $p\in[1,\infty)$. In particular,
\begin{equation}
\label{eq:rec1}
\rho(\f_\ell)=\fintl_{B^k_1}\f_\ell\,dt=1+o(1)\,,\quad \fintl_{B^k_1}|\f_\ell|^2\,dt=1+o(1)\,.
\end{equation}
By direct computation one gets
\begin{equation}\label{eq:rec1new}
\|\nabla \f_\ell\|_2^2\le c~\!\ell
\end{equation}
for $\ell\ge 1$ where, here and below, $c$ is a constant not depending on $\ell$.
Since
$$
\begin{aligned}
\mathcal E^s_k(\f_\ell)=&~\intl_{\R^k}|\tau|^{2s}|\mathcal F_k[\f_\ell]|^2~\!d\tau\le 
\Big(\intl_{\R^k}|\mathcal F_k[\f_\ell]|^2~\!d\tau\Big)^{1-s}\Big(\intl_{\R^k}|\tau|^{2}|\mathcal F_k[\f_\ell]|^2~\!d\tau\Big)^{s}\\
=&~
\|\f_\ell\|^{2(1-s)}_{L^2(B^k_1)}\|\nabla \f_\ell\|^{2s}_{L^2(B^k_1)}
\end{aligned}
$$
by the Parseval formula and by the Cauchy-Schwarz inequality, we use \eqref{eq:rec1} and \eqref{eq:rec1new} to infer that 
\begin{equation}
\label{eq:rec2}
\mathcal E^s_k(\f_\ell)\le c~\!\ell^s
\end{equation}
for $\ell\ge 1$. 

Now we prove that 
$$
\overline{v}_\ell(x,t)=~\!u(x)~\!{\f_\ell(t)}~,\quad \overline{v}_\ell\in \widetilde H^s(\Omega_1^{n+k}),
$$
is indeed a recovery sequence for $u$. First of all, 
$\rho(\overline{v}_\ell)=u~\!\rho(\f_\ell)\to u$ in $\widetilde H^s(\omega^n)$ because of the first limit in (\ref{eq:rec1}).
Secondly, from (\ref{eq:second_bis}) in Lemma \ref{L:dim} we obtain
$$
\begin{aligned}
\Iscr_\ell(\overline{v}_\ell)=&~
\frac{1}{|B^k_\ell|}\E^s_{n+k}(u\breve{\f}_\ell)\le \frac{1}{|B^k_\ell|}\widetilde{\mathcal E}^s_{n+k}(u\breve{\f}_\ell)
=
\E^s_{n}(u)\fintl_{B^k_\ell}|\breve{\f}_\ell|^2~\!dt+\frac{1}{|B^k_\ell|}\E^s_{k}(\breve{\f}_\ell)\intl_{\R^n}|u|^2~\!dx
\\
=&~
\E^s_{n}(u)\fint\limits_{B^k_1}|\f_\ell|^2~\!dt+\frac{c}{\ell^{2s}}\E^s_{k}(\f_\ell)\intl_{\R^n}|u|^2~\!dx=
\E^s_{n}(u)+o(1)~\!,
\end{aligned}
$$
where we used also (\ref{eq:rec1}) and (\ref{eq:rec2}).  We infer that
$\limsup\limits_{\ell\to \infty} \Iscr_\ell(\overline{v}_\ell)\le \E^s_{n}(u)$. Finally, noticing that $\f_\ell u\to u$ in $L^2(\Omega_1^{n+k})$
and using the assumption (\ref{eq:ip1}), we readily get
$$
\begin{aligned}
\G_\ell(\overline{v}_\ell)=&~\frac{1}{|B^k_1|}\intl_{\Omega_1^{n+k}} \invbreve{f}_\ell(x,t)\f_\ell(t)u(x)~\!dxdt
=\frac{1}{|B^k_1|}\intl_{\Omega_1^{n+k}} f_\infty(x) u(x)~\!dxdt+o(1)\\
=&~
\frac{1}{|B^k_1|}\intl_{B^k_1}dt \intl_{\omega^n}\f_\ell(x)u(x)~\!dx+o(1)=
\mathcal G_\infty(u)+o(1)~\!,
\end{aligned}
$$
and the limsup inequality in $ii)$  is proved.
This ends the proof.
\QED

In the original De Giorgi's theory,  the convergence of minima and minimizers 
of a $\Gamma$-converging sequence of equicoercive functionals is a basic result. 
Since we are not using the standard definition of $\Gamma$-convergence, we cannot relay on such general results and have to provide an {\em ad hoc} proof.
   
\bigskip

\noindent
{\bf Proof of Theorem \ref{T:Gaeta02}}
Recall that $u_\ell=\breve{v}_\ell$, where $v_\ell$ is the unique minimizer for $\invbreve{M}_\ell$. Since
$\rho_\ell(u_\ell)=\rho(v_\ell)$,  it suffices to prove that
$$
\rho(v_\ell) \to u_\infty\qquad\text{strongly in $\widetilde H^s(\omega^n)$.}
$$
Let $\overline v_\ell\in \Ht^s(\Omega_1^{n+k})$ be a {recovery sequence} for the minimizer $u_\infty$. Then
\begin{equation}
\label{eq:quasi}
\limsup_{\ell\to\infty}\invbreve{M}_\ell \le \limsup_{\ell\to\infty} I_\ell(\overline{v}_\ell)
\le I_\infty(u_\infty)~\!.
\end{equation}
Next, take any subsequence $\ell_j\to\infty$. From (\ref{eq:quasi}) we have that the sequences $\invbreve{M}_{\ell_{j}}$ and $I_{\ell_{j}}(\overline{v}_{\ell_{j}})$ are bounded.
By Lemma \ref{L:equicoercive}, 
there exists a non relabeled subsequence $v_{\ell_j}$ and $\hat v\in \Ht^s(\omega^n)$  such that $\rho(v_{\ell_j}) \to \hat v$ weakly in $\Ht^s(\omega^n)$.
By the liminf inequality in Theorem \ref{T:Gamma_omega}  we have
$$
M_\infty\le I_\infty(\hat v)\le 
\liminf_{j\to\infty} I_{\ell_j}(v_{\ell_j})= \liminf_{j\to\infty} {{\invbreve{M}_{\ell_j}}}~\!,
$$
that compared with (\ref{eq:quasi}) gives $M_\infty= I_\infty(\hat v)= {{\invbreve{M}_{\ell_j}}}+o(1)$.
It follows that  $\hat v= u_\infty$, by the uniqueness of the minimizer of the limit functional, and in fact
\begin{equation}
\label{eq:emme}
\invbreve{M}_\ell\to M_\infty~,\qquad \rho(v_\ell)\to u_\infty\quad\text{ weakly in $\widetilde H^s(\omega^n)$}
\end{equation}
as $\ell\to\infty$, because the subsequence $\ell_j$ was arbitrarily chosen. Now we use the weak lower semicontinuity of the functional $\mathcal E^s_n$,  (\ref{eq:emme}), the equalities in  (\ref{eq:G}) with $v$ replaced by $u_\infty$
and (\ref{eq:I_ine})
to obtain
$$
\begin{aligned}
\liminf_{\ell\to \infty} \frac12\mathcal E^s_n(\rho(v_\ell)) \ge &~ \frac12 \mathcal E^s_n(u_\infty) \ge M_\infty+{\mathcal G}_\infty(u_\infty)\\
=&~\lim_{\ell\to\infty}\big(\invbreve{M}_\ell+\G_\ell(v_\ell)\big)=
\lim_{\ell\to\infty} \frac12\Iscr_\ell (v_\ell) \ge
\lim_{\ell\to\infty} \frac12\mathcal E^s_n(\rho(v_\ell))
\end{aligned}
$$
that, together with the weak convergence $\rho(v_\ell)\to u_\infty$ in the uniformly convex space $\widetilde H^s(\omega^n)$,
implies that $\rho(v_\ell)\to u_\infty$ strongly in $\widetilde H^s(\omega^n)$ and 
the theorem is proved. \hfill\QED

\section{Concluding remarks}\label{sec:rem}

The aim of this section is to point out few general remarks that further clarify the choice of the convergence under which we computed the $\Gamma$-limit and the use of the averaging operators.

In particular, we study the $\Gamma$-limit of the family $\Iscr_\ell$, as $\ell\to\infty$, with respect to
other two different topologies on $\widetilde H^s(\Omega^{n+k}_1)$. We start with the standard weak topology.
Let us recall the classical notion of (sequential) $\Gamma$-convergence.

\begin{Definition}
\label{Def:Gamma}
Let $(X,\sigma)$ be a first countable topological space and let $J_h:X\to  \overline\R$
be a given sequence of functionals. The functional $J_\infty:X\to \overline\R$ is said to be the $\Gamma$-limit of the
sequence $J_h$, and we write
$\displaystyle{J_\infty=\Gamma(\sigma)\lim_{h\to \infty} J_h}$, if 
\begin{itemize}
\item[$i)$] for every $v\in  X$ and every sequence $v_h\in X$ such that 
$v_h\to v$, it holds that $\displaystyle{J_\infty(u)\le\liminf_{\ell\to\infty}J_h(v_h)}$;

\item[$ii)$] for every $v\in X$ there exists a  sequence $\overline{v}_h\in X$ such that
$\overline{v}_h\to v$ and such that $\displaystyle{J_\infty(u)\ge\limsup_{h\to\infty}J_h(\overline{v}_h)}$~\!.
\end{itemize}  
\end{Definition}

\begin{Theorem}
\label{T:Gamma_w}
Let $\sigma_{\!w}$ be the weak topology on $\widetilde H^s(\Omega^{n+k}_1)$. Then 
$$
\Gamma(\sigma_{\!w})\lim_{\ell\to\infty}\Iscr_\ell(v)=\fint\limits_{B^k_1}\mathcal E^s_n(v(~\!\cdot~\!,t))~\!dt\quad\text{for any $v\in \widetilde  H^s(\Omega^{n+k}_1)$.}
$$
\end{Theorem}

\proof
First of all we notice that the functional 
$$
\widetilde{\Iscr}_{\!\infty}(v)=\fint\limits_{B^k_1}\mathcal E^s_n(v(~\!\cdot~\!,t))~\!dt
$$
is weakly lower semicontinuous on $\widetilde  H^s(\Omega^{n+k}_1)$ because it is continuous and convex. 

We use again the transform $\breve{v}_\ell(x,t)= v\big(x,\frac{t}{\ell}\big)$ to introduce the new sequence of functionals 
$$
\widetilde{\Iscr}_{\!\ell}(v) =\frac{1}{|B^k_\ell|}~\!\widetilde\E^s_{n+k}(\breve{v}_\ell)
=\fint\limits_{B^k_1}\mathcal E^s_n(v(~\!\cdot~\!,t))~\!dt+\frac{1}{\ell^{2s}|B^k_1|}\intl_{\R^n}\mathcal E^s_k(v(x,~\!\cdot~\!))~\!dx~\!,
$$
according to \eqref{eq:tildeE}.
Trivially $\widetilde{\Iscr}_{\!\ell}(v)$ decreases to 
$\widetilde{\Iscr}_{\!\infty}(v)$
for any $v\in \widetilde  H^s(\Omega^{n+k}_1)$. By \cite{DM}, Proposition 5.7, the pointwise and the $\Gamma(\sigma_w)$-limit of the sequence
$\widetilde{\Iscr}_{\!\ell}$ coincide. 
On the other hand, from Lemma \ref{L:dim} we get
$$
\widetilde{\Iscr}_{\!\infty}(v)\le {\Iscr}_{\!\ell}(v)\le \widetilde{\Iscr}_{\!\ell}(v)\quad
\text{for any $v\in \widetilde H^s(\Omega^{n+k}_1)$,}
$$
thus the conclusion follows from \cite{DM}, Proposition 6.7.
\QED

\begin{Remark}
\label{R:Gamma_w1}
Theorem \ref{T:Gamma_w} can be equivalently stated as follows. Extend the functionals $\widetilde{\Iscr}_{\!\ell}, \widetilde{\Iscr}_{\!\infty}$
to the $L^2(\Omega_1^{n+k})$ by putting
$\widetilde{\Iscr}_{\!\ell}(v)=\widetilde{\Iscr}_{\!\infty}(v)=\infty$ if $v\notin \widetilde H^s(\Omega_1^{n+k})$ and endow
$L^2(\Omega_1^{n+k})$ with the norm topology 
$\sigma_{\!L^2}$. It holds that 
$$
\widetilde{\Iscr}_{\!\infty}= \Gamma(\sigma_{\!L^2})\lim_{\ell\to\infty}\widetilde{\Iscr}_{\!\ell}.
$$
\end{Remark}

\begin{Remark}
The limit function $\widetilde{\Iscr}_{\!\infty}$
is clearly not coercive. For instance, take
a sequence of functions $\f_\ell\in C^\infty_0(B_1)$ which is bounded in $L^2(B_1)$, but unbounded in $\widetilde H^s(B^{k}_1)$.
Then, for any nontrivial function $v\in C^\infty_0(\omega^n)$, the sequence $v\f_h$ is  unbounded  in $\widetilde H^s(\Omega_1^{n+k})$,
even if
$\widetilde{\Iscr}_{\!\infty}(v\f_h) = \mathcal E^s_n(v)~\!\|\f_h\|^2_{L^2(B_1^k)}$ is uniformly bounded.
\end{Remark}

We can partially recover coercivity by considering on $\widetilde H^s(\Omega_1^{n+k})$ a weaker topology.
For our purposes, the appropriate topology is related to the averaging operator $\rho=\rho_1$, 
$$
\rho(v)(x)=\fint\limits_{B_1^k} v(x,t)~\!dt~,\quad \rho : \Ht^s(\Omega_1^{n+k})\to \widetilde H^s(\omega^n)~\!,
$$
compare with (\ref{eq:average}) and Lemma \ref{L:average}. We denote by 
$\widetilde H^s(\omega^n)'$ the topological dual of $\widetilde H^s(\omega^n)$ and set up the next definition.

\begin{Definition}
The weak-averaged topology $\sigma_{\!\rho}$ on $\Ht^s(\Omega_1^{n+k})$ is the weakest 
 topology $\sigma$, such that
$T\circ \rho:(\Ht^s(\Omega_1^{n+k}),\sigma)\to \R $  is continuous, for any $T\in \widetilde H^s(\omega^n)'$.

We will write
$$
v_h\to_{\sigma_{\!\rho}}\!v
$$
if $v_h\to v$  in the topological space $\big(\Ht^s(\Omega_1^{n+k}),\sigma_{\!\rho}\big)$.
\end{Definition}
Clearly $\sigma_{\!\rho}$ is weaker than the 
weak topology of  $\Ht^s(\Omega_1^{n+k})$. In particular,  if $v_h \to v$ weakly in $\Ht^s(\Omega_1^{n+k})$
then $v_h\to_{\sigma_{\!\rho}}\!v$ while the converse is false, in general. 

In the next lemma we characterize convergent $\sigma_{\!\rho}$-sequences.

\begin{Lemma}
\label{L:sigma}
Let $v_h$ be a sequence in $\Ht^s(\Omega_1^{n+k})$. Then $v_h \to_{\sigma_{\!\rho}} v$ if and only if
$\rho(v_h)\to \rho(v)$ weakly in $\widetilde H^s(\omega^n)$. In this case, $\rho(v_h)\to \rho(v)$ strongly in $L^2(\omega^n)$. 
\end{Lemma}

\proof
It is well known that $v_h \to v$ in the $\sigma_{\!\rho}$-topology if and only if $(T\circ \rho)(v_h)\to (T\circ \rho)(v)$ for any 
continuous linear form $T$ on $\widetilde H^s(\omega^n)$. On the other hand, thanks to the Poincar\'e inequality on $\omega^n$, 
the immersion $\widetilde H^s(\omega^n)\to \widetilde H^s(\omega^n)'$ given by
$$
\psi\mapsto T_\psi ~,\quad \langle \Dsn \psi, v\rangle =T_\psi(v)\quad\text{for $u\in \widetilde H^s(\omega^n)$}
$$
is an invertible isometry. Finally, since
$$
(T_\psi\circ \rho)(v)= \langle \Dsn \psi, \rho(v)\rangle
$$
for any $\psi\in \widetilde H^s(\omega^n), v\in  \Ht^s(\Omega_1^{n+k})$, the first part of the statement follows. To conclude the proof,
recall that $\widetilde H^s(\omega^n)$ is compactly embedded in $L^2(\omega^n)$ by Rellich's theorem.
\QED

\medskip

We compute the 
Gamma limit of the sequence of functionals $\Iscr_\ell$ with respect to the topology $\sigma_{\!\rho}$ in the next result. The proof 
does not differ too much from that one of Theorem \ref{T:Gamma_omega}, and is omitted.

\begin{Theorem}
\label{T:Gamma_wrho}
The sequence of functionals $\Iscr_\ell:\big(\Ht^s(\Omega_1^{n+k}),\sigma_{\!\rho}\big)\to\R$, $\Gamma_{\!\sigma_{\!\rho}}$-converges as $\ell\to\infty$.
Precisely,
$$
\Gamma_{\!\sigma_{\!\rho}}\lim_{\ell\to\infty}\Iscr_\ell=~\!\mathcal E^s_n\circ \rho.
$$
\end{Theorem}

\begin{Remark}
The functional $\mathcal E^s_{n}\circ \rho: \widetilde H^s(\Omega_1^{n+k}) \to \R$ is clearly lower semicontinuous with respect to
the $\sigma_{\!\rho}$ topology, but it is not coercive. Take for instance a sequence $\ell\to \infty$ and define $\f_\ell$ as in the
proof of Theorem \ref{T:Gamma_omega}. For any nontrivial function $u\in \widetilde H^s(\omega^n)$ we have that
$\mathcal E^s_n(\rho(\f_\ell u))=\rho(\f_\ell)^2 \mathcal E^s_n(u)=\mathcal E^s_n(u)+o(1)$ is bounded; since 
$\f_\ell u\to u\notin \widetilde H^s(\Omega_1^{n+k})$ in $L^2(B^k_1)$,
then $\f_\ell u$ does not have any $\sigma_{\!\rho}$-convergent subsequence in $\widetilde H^s(\Omega_1^{n+k})$.
\end{Remark}

The last remark explains why the variant notion of $\Gamma$-convergence of \cite{ABP94} fits better to our asymptotic analysis: the restriction of the limit problem to the space $\widetilde H^s(\omega^n)$ turns out to be coercive, while its extension by $+\infty$ to the bigger space $\widetilde H^s(\Omega_1^{n+k})$ does not.

\footnotesize
\label{References}


\begin{thebibliography}{XX}
\footnotesize

\bibitem{ABP94}
G. Anzellotti, S. Baldo\ and\ D. Percivale, Dimension reduction in variational problems, 
asymptotic development in $\Gamma$-convergence and thin structures in elasticity, 
Asymptotic Anal. {\bf 9} (1994), no.~1, 61--100.

\bibitem{Braides}
A. Braides, {\it $\Gamma$-convergence for beginners}, 
Oxford Lecture Series in Mathematics and its Applications, 22, Oxford University Press, Oxford, 2002.

\bibitem{Ch}
M. Chipot, {\it $\ell$ goes to plus infinity}, 
Birkh\"{a}user Advanced Texts: Basler Lehrb\"{u}cher., Birkh\"{a}user Verlag, Basel, 2002.

\bibitem{Ch1}
M. Chipot\ and\ A. Rougirel, 
On the asymptotic behaviour of the solution of elliptic problems in cylindrical domains becoming unbounded, 
Commun. Contemp. Math. {\bf 4} (2002), no.~1, 15--44. 

\bibitem{Ch2}
M. Chipot, P. Roy\ and\ I. Shafrir, 
Asymptotics of eigenstates of elliptic problems with mixed boundary data on domains tending to infinity, 
Asymptot. Anal. {\bf 85} (2013), no.~3-4, 199--227. 

\bibitem{Ch3}
M. Chipot and \ A.Rougirel, 
On the asymptotic behaviour of the eigenmodes for elliptic problems in domains becoming unbounded, 
Trans. Amer. Math. Soc. {\bf 360} (2008), no. 7, 3579--3602.

\bibitem{ChCRS}
I. Chowdhury, G. Csat\'o, P. Roy, F. Sk,
Study of fractional Poincar\'e inequalities on unbounded domains, arXiv:1904.07170.

\bibitem{ChowR}
I. Chowdhury\ and\ P. Roy, 
On the asymptotic analysis of problems involving fractional Laplacian in cylindrical domains tending to infinity, Commun. Contemp. Math. {\bf 19} (2017), no.~5, 1650035, 21 pp.

\bibitem{CiarletII} P. G. Ciarlet, {\it Mathematical elasticity. Vol. II. Theory of plates}, Studies in Mathematics and its Applications, 27, North-Holland Publishing Co., Amsterdam, 1997.

\bibitem{DM}
G. Dal Maso, {\it An introduction to $\Gamma$-convergence}, 
 Progress in Nonlinear Differential Equations and their Applications, 8, Birkh\"{a}user Boston, Inc., Boston, MA, 1993.

\bibitem{DeGiorgi}
E. De Giorgi, $G$-convergenza e $G$-convergenza, Boll. Un. Mat. Ital. A (5) {\bf 14} (1977), no.~2, 213--220.

\bibitem{FP04} L. Freddi\ and\ R. Paroni, 
The energy density of martensitic thin films via dimension reduction, 
Interfaces Free Bound. {\bf 6} (2004), no.~4, 439--459.

\bibitem{FP04bis} L. Freddi\ and\ R. Paroni, 
A 3D--1D Young measure theory of an elastic string, Asymptot. Anal. {\bf 39} (2004), no.~1, 61--89. 


\bibitem{K}
M. Kwa\'{s}nicki, Ten equivalent definitions of the fractional Laplace operator, Fract. Calc. Appl. Anal. {\bf 20} (2017), no.~1, 7--51. 

\bibitem{Love}
A. E. H. Love, {\it A treatise on the Mathematical Theory of Elasticity}, 
Dover Publications, New York, 1944.

\bibitem{MN}
R. Musina\ and\ A. I. Nazarov, 
On the Sobolev and Hardy constants for the fractional Navier Laplacian, 
Nonlinear Anal. {\bf 121} (2015), 123--129.

\bibitem{MaxPri}
R. Musina\ and\ A. I. Nazarov, 
Strong maximum principles for fractional Laplacians, 
Proc. Roy. Soc. Edinburgh Sect. A {\bf 149} (2019), no.~5, 1223--1240.

\bibitem{SV}
A. J. C. B. Saint-Venant, Memoire sur la Torsion des Prismes, Mem. Divers Savants {\bf 14} (1855), 233--560.

\bibitem{SeVa}
R. Servadei\ and\ E. Valdinoci, 
Variational methods for non-local operators of elliptic type, 
Discrete Contin. Dyn. Syst. {\bf 33} (2013), no.~5, 2105--2137.

\bibitem{Sil}
L. Silvestre, Regularity of the obstacle problem for a fractional power of the Laplace operator, 
Comm. Pure Appl. Math. {\bf 60} (2007), no.~1, 67--112.

\bibitem{Toupin}
R. A. Toupin, Saint-Venant's principle, Arch. Rational Mech. Anal. {\bf 18} (1965), 83--96.

\bibitem{Tr}
H. Triebel, {\it Interpolation theory, function spaces, differential operators}, 
VEB Deutscher Verlag der Wissenschaften, Berlin, 1978.

\bibitem{Y}
K. Yeressian, 
K. Yeressian, Asymptotic behavior of elliptic nonlocal equations set in cylinders,
Asymptot. Anal. {\bf 89} (2014), no.~1-2, 21--35. 

\end{thebibliography}
\end{document}